\documentstyle[11pt,amssymb,amsfonts,twoside]{article}

\begin{titlepage}

\title{Index and secondary index theory for flat bundles with duality}
\author{Ulrich Bunke\thanks{G\"ottingen, bunke@uni-math.gwdg.de} and Xiaonan
Ma\thanks{Berlin, ma@mathematik.hu-berlin.de}}

\end{titlepage}

\begin{document}

\maketitle
\newcommand{\diag}{{\rm diag}}
\newcommand{\proof}{{\it Proof.$\:\:\:\:$}}
 \newcommand{\dist}{{\rm dist}}
\newcommand{\kaaa}{{\frak k}}
\newcommand{\paaa}{{\frak p}}
\newcommand{\vp}{{\varphi}}
\newcommand{\taaa}{{\frak t}}
\newcommand{\haaa}{{\frak h}}
\newcommand{\R}{{\Bbb R}}
\newcommand{\Hh}{{\bf H}}
\newcommand{\Rep}{{\rm Rep}}
\newcommand{\Hb}{{\Bbb H}}
\newcommand{\Q}{{\Bbb Q}}
\newcommand{\str}{{\rm str}}
\newcommand{\Ind}{{\rm ind}}
\newcommand{\triv}{{\rm triv}}
\newcommand{\Z}{{\Bbb Z}}
\newcommand{\bD}{{\bf D}}
\newcommand{\bF}{{\bf F}}
\newcommand{\tX}{{\tt X}}
\newcommand{\Cliff}{{\rm Cliff}}
\newcommand{\tY}{{\tt Y}}
\newcommand{\tZ}{{\tt Z}}
\newcommand{\tV}{{\tt V}}
\newcommand{\tR}{{\tt R}}
\newcommand{\Fam}{{\rm Fam}}
\newcommand{\Cusp}{{\rm Cusp}}
\newcommand{\bT}{{\bf T}}
\newcommand{\bK}{{\bf K}}
\newcommand{\K}{{\Bbb K}}
\newcommand{\tH}{{\tt H}}
\newcommand{\bS}{{\bf S}}
\newcommand{\bB}{{\bf B}}
\newcommand{\tW}{{\tt W}}
\newcommand{\tF}{{\tt F}}
\newcommand{\bA}{{\bf A}}
\newcommand{\bL}{{\bf L}}
 \newcommand{\bom}{{\bf \Omega}}
\newcommand{\bundle}{{bundle}}
\newcommand{\ch}{{\bf ch}}
\newcommand{\ve}{{\varepsilon}}
\newcommand{\C}{{\Bbb C}}
\newcommand{\gen}{{\rm gen}}
\newcommand{\cTop}{{{\cal T}op}}
\newcommand{\bP}{{\bf P}}
\newcommand{\Naaa}{{\bf N}}
\newcommand{\image}{{\rm image}}
\newcommand{\gaaa}{{\frak g}}
\newcommand{\zaaa}{{\frak z}}
\newcommand{\saaa}{{\frak s}}
\newcommand{\laaa}{{\frak l}}
\newcommand{\stimes}{{\times\hspace{-1mm}\bf |}}
\newcommand{\ausg}{{\rm end}}
\newcommand{\bff}{{\bf f}}
\newcommand{\maaa}{{\frak m}}
\newcommand{\aaaa}{{\frak a}}
\newcommand{\naaa}{{\frak n}}
\newcommand{\brr}{{\bf r}}
\newcommand{\res}{{\rm res}}
\newcommand{\Aut}{{\rm Aut}}
\newcommand{\Pol}{{\rm Pol}}
\newcommand{\Tr}{{\rm Tr}}
\newcommand{\cT}{{\cal T}}
\newcommand{\dom}{{\rm dom}}
\newcommand{\db}{{\bar{\partial}}}
\newcommand{\g}{{\gaaa}}
\newcommand{\cZ}{{\cal Z}}
\newcommand{\cH}{{\cal H}}
\newcommand{\cM}{{\cal M}}
\newcommand{\interi}{{\rm int}}
\newcommand{\singsupp}{{\rm singsupp}}
\newcommand{\cE}{{\cal E}}
\newcommand{\ccR}{{\cal R}}
\newcommand{\cV}{{\cal V}}
\newcommand{\cY}{{\cal Y}}
\newcommand{\cW}{{\cal W}}
\newcommand{\cI}{{\cal I}}
\newcommand{\cC}{{\cal C}}
\newcommand{\mod}{{\rm mod}}
\newcommand{\cK}{{\cal K}}
\newcommand{\cA}{{\cal A}}
\newcommand{\cEp}{{{\cal E}^\prime}}
\newcommand{\cU}{{\cal U}}
\newcommand{\Hom}{{\mbox{\rm Hom}}}
\newcommand{\vol}{{\rm vol}}
\newcommand{\cO}{{\cal O}}
\newcommand{\End}{{\mbox{\rm End}}}
\newcommand{\Ext}{{\mbox{\rm Ext}}}
\newcommand{\rk}{{\rm rank}}
\newcommand{\im}{{\mbox{\rm im}}}
\newcommand{\sign}{{\rm sign}}
\newcommand{\spann}{{\mbox{\rm span}}}
\newcommand{\symm}{{\mbox{\rm symm}}}
\newcommand{\cF}{{\cal F}}
\newcommand{\cD}{{\cal D}}
\newcommand{\Ree}{{\rm Re }}
\newcommand{\Res}{{\mbox{\rm Res}}}
\newcommand{\Imm}{{\rm Im}}
\newcommand{\inter}{{\rm int}}
\newcommand{\clo}{{\rm clo}}
\newcommand{\tg}{{\rm tg}}
\newcommand{\ee}{{\rm e}}
\newcommand{\Li}{{\rm Li}}
\newcommand{\cN}{{\cal N}}
 \newcommand{\conv}{{\rm conv}}
\newcommand{\op}{{\mbox{\rm Op}}}
\newcommand{\tr}{{\mbox{\rm tr}}}
\newcommand{\cs}{{c_\sigma}}
\newcommand{\ctg}{{\rm ctg}}
\newcommand{\degg}{{\mbox{\rm deg}}}
\newcommand{\Ad}{{\mbox{\rm Ad}}}
\newcommand{\ad}{{\mbox{\rm ad}}}
\newcommand{\codim}{{\rm codim}}
\newcommand{\Gr}{{\mathrm{Gr}}}
\newcommand{\coker}{{\rm coker}}
\newcommand{\id}{{\mbox{\rm id}}}
\newcommand{\ord}{{\rm ord}}
\newcommand{\nat}{{\Bbb  N}}
\newcommand{\supp}{{\rm supp}}
\newcommand{\sing}{{\mbox{\rm sing}}}
\newcommand{\spec}{{\mbox{\rm spec}}}
\newcommand{\Ann}{{\mbox{\rm Ann}}}
\newcommand{\aca}{{\aaaa_\C^\ast}}
\newcommand{\acag}{{\aaaa_{\C,good}^\ast}}
\newcommand{\acage}{{\aaaa_{\C,good}^{\ast,extended}}}
\newcommand{\tck}{{\tilde{\ck}}}
\newcommand{\tnk}{{\tilde{\ck}_0}}
\newcommand{\ceep}{{{\cal E}(E)^\prime}}
 \newcommand{\ncE}{{{}^\naaa\cE}}
 \newcommand{\Or}{{\rm Or }}
\newcommand{\Diff}{{\cal D}iff}
\newcommand{\cB}{{\cal B}}
\newcommand{\hc}{{{\cal HC}(\gaaa,K)}}
\newcommand{\hcma}{{{\cal HC}(\maaa_P\oplus\aaaa_P,K_P)}}
\def\imath{{\rm i}}
\newcommand{\vsl}{{V_{\sigma_\lambda}}}
\newcommand{\czg}{{\cZ(\gaaa)}}
\newcommand{\csl}{{\chi_{\sigma,\lambda}}}
\newcommand{\cR}{{R}}
\def\hB{\hspace*{\fill}$\Box$ \newline\noindent}
\newcommand{\varho}{\varrho}
\newcommand{\ind}{{\rm index}}
\newcommand{\Indu}{{\rm Ind}}
\newcommand{\Fin}{{\mbox{\rm Fin}}}
\newcommand{\cS}{{S}}
\newcommand{\orig}{{\cal O}}
\def\hB{\hspace*{\fill}$\Box$ \\[0.5cm]\noindent}
\newcommand{\cL}{{\cal L}}
 \newcommand{\cG}{{\cal G}}
\newcommand{\npci}{{\naaa_P\hspace{-1.5mm}-\hspace{-2mm}\mbox{\rm coinv}}}
\newcommand{\pki}{{(\paaa,K_P)\hspace{-1.5mm}-\hspace{-2mm}\mbox{\rm inv}}}
\newcommand{\mki}{{(\maaa_P\oplus \aaaa_P, K_P)\hspace{-1.5mm}-\hspace{-2mm}\mbox{\rm inv}}}
\newcommand{\Mat}{{\rm Mat}}
\newcommand{\npi}{{\naaa_P\hspace{-1.5mm}-\hspace{-2mm}\mbox{\rm inv}}}
\newcommand{\ngp}{{N_\Gamma(\pi)}}
\newcommand{\gbg}{{\Gamma\backslash G}}
\newcommand{\gkm}{{ Mod(\gaaa,K) }}
\newcommand{\ggkm}{{  (\gaaa,K) }}
\newcommand{\pkm}{{ Mod(\paaa,K_P)}}
\newcommand{\ppkm}{{  (\paaa,K_P)}}
\newcommand{\makm}{{Mod(\maaa_P\oplus\aaaa_P,K_P)}}
\newcommand{\mmakm}{{ (\maaa_P\oplus\aaaa_P,K_P)}}
\newcommand{\cP}{{\cal P}}
\newcommand{\gm}{{Mod(G)}}
\newcommand{\gk}{{\Gamma_K}}
\newcommand{\La}{{\cal L}}
\newcommand{\cug}{{\cU(\gaaa)}}
\newcommand{\cuk}{{\cU(\kaaa)}}
\newcommand{\dc}{{C^{-\infty}_c(G) }}
\newcommand{\gdk}{{\gaaa/\kaaa}}
\newcommand{\dgkm}{{ D^+(\gaaa,K)-\mbox{\rm mod}}}
\newcommand{\dgm}{{D^+G-\mbox{\rm mod}}}
\newcommand{\vect}{{\C-\mbox{\rm vect}}}
 \newcommand{\cig}{{C^{ \infty}(G)_{K} }}
\newcommand{\gami}{{\Gamma\hspace{-1.5mm}-\hspace{-2mm}\mbox{\rm inv}}}
\newcommand{\cQ}{{\cal Q}}
\newcommand{\mmap}{{Mod(M_PA_P)}}
\newcommand{\bbbz}{{\bf Z}}
 \newcommand{\cX}{{\cal X}}
\newcommand{\bH}{{\bf H}}
\newcommand{\pr}{{\rm pr}}
\newcommand{\bX}{{\bf X}}
\newcommand{\bY}{{\bf Y}}
\newcommand{\bZ}{{\bf Z}}
\newcommand{\bV}{{\bf V}}

\newtheorem{prop}{Proposition}[section]
\newtheorem{lem}[prop]{Lemma}
\newtheorem{ddd}[prop]{Definition}
\newtheorem{theorem}[prop]{Theorem}
\newtheorem{kor}[prop]{Corollary}
\newtheorem{ass}[prop]{Assumption}
\newtheorem{con}[prop]{Conjecture}
\newtheorem{prob}[prop]{Problem}
\newtheorem{fact}[prop]{Fact}
\begin{abstract}
We discuss some aspects of index and secondary index theory for flat bundles
with duality. This theory was first developed by Lott \cite{lott99}. Our
main purpose in the present paper is provide a modification with better
functorial properties. \end{abstract}

\tableofcontents

\section{The functor $L$}

\subsection{Definition of $L$}

We define a contravariant functor $L$ from the category $\cTop$ of topological
spaces and continuous maps to $\Z_2$-graded rings. 

Let $X$ be a topological space and $\epsilon\in \Z_2$.
If $\cF$ is a locally constant sheaf of finite-dimensional $\R$-modules over
$X$, then an {\bf $\epsilon$-symmetric duality structure} on $\cF$ is an
isomorphism of sheafs $q:\cF\rightarrow \cF^*$ satisfying $q^*=(-1)^\epsilon
q$, where $\cF^*:=\Hom(\cF,\underline{\R})$.

We first define a semigroup $\hat L_\epsilon(X)$ and then obtain
$L_\epsilon(X)$ by introducing a relation. A generator of the semigroup
$\hat L_\epsilon(X)$ is a pair $(\cF,q)$ consisting of a locally constant sheaf
of finite-dimensional $\R$-modules and an $\epsilon$-symmetric duality
structure $q$. The operation in $\hat L_\epsilon(X)$ is given by direct sum
 $$(\cF,q)+(\cF^\prime,q^\prime):=(\cF\oplus \cF^\prime,q\oplus q^\prime)\
.$$ 
The group $L_\epsilon(X)$ is obtained from $\hat L_\epsilon(X)$ by introducing
the  relation {\bf lagrangian reduction}. 
We require $(\cF,q)\sim 0$ if there exists a  locally constant lagrangian
subsheaf $i:\cL\hookrightarrow  \cF$, i.e.  the composition
$\cL\stackrel{i}{\rightarrow} \cF\stackrel{q}{\rightarrow} \cF^*
\stackrel{i^*}{\rightarrow} \cL^*$ vanishes and $\cL=\cL^\perp:=\ker(i^*\circ
q)$.  

The class of $(\cF,q)$ in $L_\epsilon(X)$ will be denoted by $[\cF,q]$.
The ring structure $L_\epsilon(X)\otimes L_{\epsilon^\prime}(X)\rightarrow
L_{\epsilon+\epsilon^\prime}(X)$ is induced by the tensor product:
$$[\cF,q][\cF^\prime,q\prime]:=
\frac{\sqrt{(-1)^\epsilon}\sqrt{(-1)^{\epsilon^\prime}}}{\sqrt{(-1)^{\epsilon+\epsilon^\prime}} }
[\cF\otimes \cF^\prime,q\otimes q^\prime]\ .$$  
The sign-convention is made such that later we have a natural transformation
of rings from $L$ to complex $K$-theory $K^0$.

If $f:Y\rightarrow X$ is a morphism in $\cTop$, then $f^*:L(X)\rightarrow
L(Y)$ is defined by $f^*[\cF,q]=[f^*\cF, f^*q]$. The map $f^*$ only depends on
the homotopy class of $f$.

\underline{Remark:} A version $L(X)^{Lott}$ of this ring was first introduced
by Lott \cite{lott99}. His definition differs from ours since our relation
"lagrangian reduction" is repaced by "hyperbolic reduction" in the definition
of Lott. Here a generator $(\cF,q)$ is called hyperbolic if there is an
lagrangian subsheaf $\cL\subset \cF$ which
admits a complement inside $\cF$. In particular, $L(X)$ is a quotient of
$L(X)^{Lott}$.

\subsection{Computation of $L(X)$}

Assume that $X$ is path connected and fix a base point $x$. Let $G:=\pi_1(X,x)$
be the fundamental group. We consider the semigroup $\hat L_\epsilon(G)$ of
all tuples $(F,q,\rho)$, where $F$ is a finite-dimensional real vector space,
$q:F\rightarrow F^*$ is an $\epsilon$-symmetric duality structure, and
$\rho:G\rightarrow \Aut(F,q)$ is a representation of $G$. The
semigroup operation is given by direct sum. We then define the group
$L_\epsilon(G)$ by introducing the relation $(F,q,\rho)\sim 0$ if there exists
an invariant  lagrangian  subspace $L\subset F$. Taking the holonomy
representation on the fibre over $x$ we obtain a bijection $L(X)\cong L(G)$.

Let $\hat G:=G\times\Z_2$. If $(F,q,\rho)\in \hat L(G)$, then we form the
representation $(\hat F,\hat \rho)$ of $\hat G$ by 
$\hat F:=F\oplus F^*$, $\hat \rho((g,0)):=\rho(g)\oplus\rho(g^{-1})^*$,
$$\hat \rho(1,1):=\left(\begin{array}{cc}0&q^{-1}\\q&0\end{array}\right)\ .$$

If $(F,q,\rho)$ is an irreducible representation of
$G$, then $(\hat F,\hat\rho)$ is an irreducible representation of $\hat G$.
An irreducible representation $(V,\sigma)$  of $\hat G$ 
is called real, complex, or quaternionic, if
$\End_G(V)\cong \R$, $\cong \C$ or $\cong\Hb$, respectively. 

For each irreducible representation $(F,\rho)$ of $G$
we define a $\Z_2$-graded group $A(F,\rho)$ as follows.
If $(F,\rho)$ admits an invariant $\epsilon$-symmetric form $q$, then we define
$A_\epsilon(F,\rho):=\Z$, iff $(\hat F,\hat \rho)$ is real or quaternionic,
$A_\epsilon(F,\rho):=\Z_2$, iff $(\hat F,\hat \rho)$ is complex, and
we set $A_\epsilon(F,\rho):=0$, if $(F,\rho)$ does not admit an invariant 
$\epsilon$-symmetric form. If $A(F,\rho)\not=0$, then we fix an invariant
$\epsilon$-symmetric form $q_{F,\rho}$.

Let $\Rep(G)$ denote the set of isomorphism classes of irreducible
representations of $G$.  
\begin{theorem}
There is a natural  isomorphism of $\Z_2$-graded groups 
$$L(G)=\bigoplus_{(F,\rho)\in\Rep(G)} A(F,\rho)$$
fixed by the condition that $(F,q,q_{F,q})\mapsto 1_{(F,q)}$.
\end{theorem}

\subsection{The natural transformation to $K$-theory}

By $\cTop_{para}$ we denote the full subcategory of $\cTop$ of paracompact
topological spaces.  Let $K^0(X)$ be the complex $K$-theory functor. We
construct a natural transformation  $b:L\rightarrow K^0$ of functors from
$\cTop_{para}$ to rings.

A locally constant sheaf of finite-dimensional $\R$-modules on
$X$ gives rise to a locally trivial real vector bundle $\bundle(\cF)$ in a
natural way.  The correspondence $\bundle$ is functorial and compatible with
direct sum, tensor product, and duality.  Thus applying the bundle
construction to $(\cF,q)$ we obtain a pair $(F,Q)$ consisting of a
finite-dimensional real vector bundle and an isomorphism $Q:F\rightarrow F^*$.

Let $(F,Q)$ be a real vector bundle with an isomorphism $Q:F\rightarrow F^*$   
such that $Q^*=(-1)^\epsilon Q$ for $\epsilon\in \Z_2$.
Following the language introduced by Lott \cite{lott99} we call an isomorphism
$J:F\rightarrow F$ a metric structure, if
\begin{enumerate}
\item  $J^*\circ Q$ defines a scalar product on $F$,
\item $J^2=(-1)^\epsilon \id_F$,
\item $J^*\circ Q=(-1)^\epsilon Q\circ J$.
\end{enumerate}

Since we assume that $X$ is paracompact it admits partitions of unity.
This implies that  metric structures  exist, and that the space of all
metric structures is contractible. 

Given $(F,Q)$ as above we choose a metric structure $J$.
Let $F_\C$ be the complexification of $F$.
Then $\frac{1}{\sqrt{(-1)^\epsilon}} J$ is a $\Z_2$-grading of $F_\C$, and the
pair $(F_\C,\frac{1}{\sqrt{(-1)^\epsilon}} J)$ represents an element of
$K^0(X)$ which does not depend on the choice of $J$.

The transformation $b:L\rightarrow K^0$ is obtained by composing the latter
construction with $\bundle$. 

\subsection{Secondary $L$-groups. The $\R/\Z$-variant}

We first recall the definition of the $2$-periodic cohomology theory
$K_{\R/\Z}$ introduced by \cite{karoubi87}, \cite{lott9?}. Let $BU$
be the classifying space of complex $K$-theory. The Chern character
(with real coefficients)
is induced by a map $\ch:BU\rightarrow \prod_{n=1}^\infty K(\R,2n)$.
The homotopy fibre of this map classifies $K_{\R/\Z}$. In particular, for any
paracompact space $X$ there is a natural exact sequence of $K^0(X)$-modules
$$\rightarrow K^{-1}(X)\stackrel{\ch}{\rightarrow}
H^{odd}(X,\R) \rightarrow K^{-1}_{\R/\Z}(X)
\stackrel{\beta}{\rightarrow} K^0(X)\stackrel{\ch}{\rightarrow}
H^{ev}(X,\R)\rightarrow \ ,$$
where $K^0(X)$ acts on cohomology via the Chern character.

We now define the functor $X\mapsto \bar L^{\R/\Z}(X)$ from paracompact
topological spaces to $\Z_2$-graded groups by the following pull-back
diagram:   
$$\begin{array}{ccc}
\bar L^{\R/\Z}(X) & \rightarrow& L(X)\\
\downarrow&&\downarrow b\\
K^{-1}_{\R/\Z}(X)&\stackrel{\beta}{\rightarrow} & K^0(X)
\end{array}\ .
$$
The grading is induced from that of $L(X)$.
On morphisms the functor $L^{\R/\Z}$ only depends on homotopy classes.
Note that $K^0(X)$ and $K^{-1}_{\R/\Z}(X)$ are $L(X)$-modules via $b$. 
This induces a graded $L(X)$-module structure on $\bar L^{\R/\Z}(X)$.
We have the following natural exact sequence of $L(X)$-modules
$$ 
K^{-1}(X)\stackrel{\ch}{\rightarrow} H^{odd}(X,\R) \rightarrow \bar
L^{\R/\Z}(X)    \rightarrow L(X)\stackrel{\ch\circ
b}{\rightarrow} H^{ev}(X,\R)\ .$$

\subsection{Index and secondary index}

Let $X\rightarrow B$ be a smooth locally trivial fibre bundle
over a compact base $B$ with compact even-dimensional fibres such that the
vertical bundle $TX/B$ is orientable. In this case we have the following maps:
\begin{itemize}
\item 
$\pi^{sign}_*:H^*(X,\R)\rightarrow H^*(B,\R)$   defined by
$\pi_*(\omega)=\int_{X/B} \omega\cap \bL(TX/B)$,
where $\int_{X/B}$ is integration over the fibre and $\bL(TX/B)$ denotes the
Hirzebruch $\bL$-class of the vertical bundle.
\item
$\pi^{sign}_!:K^0(X)\rightarrow K^0(B)$  defined by
$\pi^{sign}([E])=\ind(D^{sign}_E)$, where $D^{sign}_E$ is the fibrewise
signature operator twisted by $E$ and $\ind(D^{sign}_E)\in K^0(B)$ denotes the
class of the index bundle.
\item $\pi^{L}_*:L(X)\rightarrow L(B)$ defined by
$\pi_*^L[\cF,q]=[H^\bullet R\pi_*\cF,\pi_*(q)]$, where
$H^\bullet R\pi_* \cF=\oplus_{i=0}^\infty R^i \pi_*(\cF)$ is the direct sum of
higher derived direct images of $\cF$. The orientation of $TX/B$ provides an
isomorphism $p:R^{\dim(TX/B)}\pi_* \underline{\R}\stackrel{\sim}{\rightarrow}
\underline{\R}$. It induces isomorphisms  
$R^i\pi_* \cF^* \rightarrow (R^{\dim(TX/B)-i}\pi_*\cF)^*$,
and  $\pi_*(q)$  is the composition of the sum of these maps
with the sum of $R^i\pi_*(q):R^i\pi_*\cF\rightarrow R^i \pi_* \cF^*$.  
\end{itemize}
It is an excercise in spectral sequences to show that $\pi^L_*$ is
well-defined. 

\underline{Remark :}
In \cite{lott99} Lott defines $\pi^{L,Lott}_*:L^{Lott}(X)\rightarrow
L^{Lott}(B)$ which induces $\pi^L_*$ by passing to  quotients.

By the index theorem for families and fibrewise Hodge theory
the following diagram commutes
$$\begin{array}{ccccc}
L(X)&\rightarrow &K^0(X)&\rightarrow& H^{ev}(X,\R)\\
\pi_*^L \:\downarrow&&\pi^{sign}_!\downarrow &&\pi^{sign}_* \:\downarrow\\
L(B)&\rightarrow &K^0(B)&\rightarrow& H^{ev}(B,\R)
\end{array}\ .$$
In order to define an index map for $K_{\R/\Z}$ we need the further assumption
that $\pi$ is $K$-oriented. Thus assume that $TX/B$ has a $Spin_c$-structure.
Then there are maps $\pi^{Spin_c}_!:K^0(X)\rightarrow K^0(B)$ and
$\pi^{Spin_c,\R/\Z}_!:K^{-1}_{\R/\Z}(X)\rightarrow K^{-1}_{\R/\Z}(B)$
(compare e.g. \cite{lott9?}), such that the following diagram commutes:
$$\begin{array}{ccccccc}
H^{odd}(X,\R)&\rightarrow&K^{-1}_{\R/\Z}(X)&\rightarrow &K^0(X)&\rightarrow&
H^{ev}(X,\R)\\ \pi^{Spin_c}_*\downarrow&&
\pi_!^{Spin_c,\R/\Z}\:\downarrow&&\pi^{Spin_c}_!\downarrow && \pi^{Spin_c}_*
\:\downarrow\\ H^{odd}(B,\R)&\rightarrow &K_{\R/\Z}^{-1}(B)&\rightarrow
&K^0(B)&\rightarrow& H^{ev}(B,\R) \end{array}\ ,$$
where $\pi^{Spin_c}_*(\omega)=\int_{X/B} \hat \bA(TX/B)\cap \ee^{c_1/2}\cap
\omega$ and $c_1$ is the first Chern class determined by the
$Spin_c$-structure.
There is an unique element $E_{sign}\in K_0(X)$ such
that $\pi_!^{sign}(x)=\pi^{Spin_c}_!(E_{sign} \bullet x)$. Note that
$\ch(E_{sign})\cap \hat \bA(TX/B)\cap \ee^{c_1/2}=\bL(TX/B)$.

\begin{itemize}
\item
$\pi^{sign,\R/\Z}:K^{-1}_{\R/\Z}(X)\rightarrow K^{-1}_{\R/\Z}(B)$ is defined by
$\pi^{sign,\R/\Z}(x):=\pi^{Spin_c,\R/\Z}(E_{sign}\bullet x)$ so
that the following diagram commutes
$$\begin{array}{ccccccc}
H^{odd}(X,\R)&\rightarrow&K^{-1}_{\R/\Z}(X)&\rightarrow &K^0(X)&\rightarrow&
H^{ev}(X,\R)\\ \pi^{sign}_*\downarrow&&
\pi_!^{sign,\R/\Z}\:\downarrow&&\pi^{sign}_!\downarrow && \pi^{sign}_*
\:\downarrow\\ H^{odd}(B,\R)&\rightarrow &K_{\R/\Z}^{-1}(B)&\rightarrow
&K^0(B)&\rightarrow& H^{ev}(B,\R) \end{array}\ .$$
\item 
$\pi^{\bar L,\R/\Z}_*:\bar L^{\R/\Z}(X)\rightarrow \bar L^{\R/\Z}(B)$ is the
map induced by $\pi^{sign,\R/\Z}_!$ and $\pi^{L}_*$. The following diagramm
commutes:
$$\begin{array}{ccccccc}
H^{odd}(X,\R)&\rightarrow&\bar L^{\R/\Z}(X)&\rightarrow &L(X)&\rightarrow&
H^{ev}(X,\R)\\ \pi^{sign}_*\downarrow&&
\pi_*^{\bar L,\R/\Z}\:\downarrow&&\pi^{L}_*\downarrow && \pi^{sign}_*
\:\downarrow\\ H^{odd}(B,\R)&\rightarrow &\bar L^{\R/\Z}(B)&\rightarrow
&L(B)&\rightarrow& H^{ev}(B,\R) \end{array}\ .$$
\end{itemize}

All index maps are natural with respect to pull-back of fibre bundles.

Let $\pi_1:X\rightarrow X_1$ and $\pi_2:X_1\rightarrow B$ be locally
trivial smooth fibre bundles with closed even-dimensional fibres and compact
base. Further assume that the vertical bundles $TX_1/X_2$ and $TX_2/B$ carry
$Spin_c$-structures (and are therefore oriented).
Then the composition $\pi=\pi_2\circ\pi_1:X\rightarrow B$ is a locally trivial
fibre bundle with closed even-dimensional fibres, and the vertical bundle
$TX/B$ carries an induced $Spin_c$-structure. In this situation the 
 index maps on complex $K$-theory, $K_{\R/\Z}$-theory,  and in  cohomology
behave functorially with respect to the iterated fibre bundle. It is again an
excercise in spectral sequences to show that $\pi^L_*$ is also functorial. As a
consequence, $\pi^{\bar L,\R/\Z}_*$ is functorial, i.e. $\pi^{\bar
L,\R/\Z}_*=(\pi_2^{\bar L,\R/\Z})_*\circ (\pi_1^{\bar L,\R/\Z})_*$.

\underline{Remark :}
The motivation of introducing the quotient $L(X)$ of Lott's group
$L^{Lott}(X)$ is to implement this functoriality for $\pi^L_*$, which is not
true for Lott's definition.

\section{Geometry and extended secondary $L$-groups}\label{geo}

\subsection{Definition of $\bar L$}

The functor $X\mapsto \bar L^{\R/\Z}(X)$ from $\cTop_{para}$
to $\Z_2$-graded $L(X)$-modules was defined by a purely
homotopy-theoretic construction as an extension of the functor $$X\mapsto
\ker\left(\ch\circ b:L(X)\rightarrow H^{ev}(X,\R)\right)$$ by  $X\mapsto 
H^{odd}(X,\R)/\ch(K^{-1}(X))$.

Let  $\cTop_{smooth}$ denote
 the full subcategory of topological spaces $\cTop_{para}$ which
are homotopy equivalent to smooth manifolds.
 In the present section we use a
differential geometric construction in order to define on $\cTop_{smooth}$ a
functor $X\mapsto \bar L(X)$ to graded $L(X)$-modules 
which extends $X\mapsto \ker
(\ch\circ b)$ by $X\mapsto  H^{odd}(X,\R)$.      
First we define $\bar L(M)$ for a smooth manifold $M$.  
Again we first define a semigroup $\hat{ \bar L}(M)$.
A generator of $\hat{\bar L}_\epsilon(M)$ is a tuple $(\cF,q,J,\rho)$, where
$(\cF,q)$ is a generator of $\hat L_\epsilon(M)$, $J$ is a metric structure on
$F:=\bundle(\cF)$, and $\rho\in \Omega^{4*-(-1)^\epsilon}(M)/\im(d)$
satisfies $d\rho=p(\nabla^F,J)$.
Here $(\Omega^*(M),d)$ is the  real deRham complex of $M$ and  $\nabla^F$ is
the canonical flat connection on $F$ (such that $\ker(\nabla^F)\cong \cF$).
The metric structure $J$ induces a $\Z_2$-grading  of $F_\C$, and
if $\nabla^{F_\C,J}$ denotes the even part of the extension of $\nabla^F$
to $F_\C$, then $p(\nabla^F,J):=\ch(\nabla^{F_\C,J})=\tr_s
\exp(-R^{\nabla^{F_\C,J}}/2\pi\imath)$ is the characteristic form
representing $\ch\circ b([\cF,q])$. The semigroup operation is given by
$$(\cF,q,J,\rho)+
(\cF^\prime,q^\prime,J^\prime,\rho^\prime):= (\cF\oplus \cF^\prime,q\oplus
q^\prime,J\oplus J^\prime,\rho+\rho^\prime)\ .$$ We form $\bar L_\epsilon(M)$
by introducing the relation {\bf lagrangian reduction}.
Let $\cL\subset \cF$ be a locally constant lagrangian subsheaf and
$L:=\bundle(\cL)$.  Then we have a decomposition $F=L\oplus  J(L)$.
Let $\nabla^\oplus$ be the part of $\nabla^{F_\C,J}$ which preserves this
decomposition, and let $\tilde\ch(\nabla^\oplus,\nabla^{F_\C,J})$ be the
transgression Chern form such that
$d\tilde\ch(\nabla^\oplus,\nabla^{F_\C})=\ch(\nabla^\oplus)-\ch(\nabla^{F_\C})$.
Note that $\ch(\nabla^\oplus)=0$.
In $\bar L_\epsilon(M)$ we require
$$(\cF,q,J,\rho)\sim  (0,0,0,\rho+\tilde \ch(\nabla^\oplus,\nabla^{F_\C,J})) \
.$$ By
$[\cF,q,J,\rho]$ we denote the  class in $\bar L$ represented by
$(\cF,q,J,\rho)$.

\underline{Remark :} In \cite{lott99} similar functors $\bar L^{Lott}$
were defined replacing the relation "lagrangian reduction" by "hyperbolic
reduction" and "change of metric structure".  Lott's relation is smaller
than ours, and we have a natural surjective map $\bar L^{Lott}(M)\rightarrow
\bar L(M)$.

The graded module structure of $\bar L(X)$ over $L(X)$ is defined by
$$[\cF,q,J,\rho]\bullet
[\cE,p]:=\frac{\sqrt{(-1)^\epsilon}\sqrt{(-1)^{\epsilon^\prime}}}{\sqrt{(-1)^{\epsilon+\epsilon^\prime}}}[\cF\otimes
\cE,q\otimes p,J\otimes J^E,\rho\wedge \ch(\nabla^{E_\C,J^E})]\ ,$$ where $J^E$
is any metric structure on $\bundle(\cE)$.

\subsection{Functorial properties}

If $f:M\rightarrow N$ is a smooth map of manifolds, then we obtain an induced
map $f^*:\bar L(N)\rightarrow \bar L(M)$, which is given on generators by
pull-back of structures.  $f^*$ only depends
on the smooth homotopy class of $f$ and is compatible with the $L(M)$-module structures.

We define natural maps $H^{odd}(M,\R)\rightarrow \bar L(M)$ and
$\bar L(M)\rightarrow L(M)$ by
$[\rho]\mapsto [0,0,0,\rho]$ and $[\cF,q,J,\rho]\mapsto [\cF,q]$.
Then we have the following exact sequence of $L(M)$-modules
$$H^{odd}(M,\R)\rightarrow \bar L(M)\rightarrow
L(M)\rightarrow H^{ev}(M,\R)\ $$
(see \cite{lott99}, Prop. 21, for a similar argument).

\begin{prop}
The map $H^{odd}(M,\R)\rightarrow \bar L(M)$ is injective.
\end{prop}
Indeed, 
let $\omega\in \Omega^{odd}(M)$ be a closed form.
If $(0,0,0,\omega)\sim 0$ in $\bar L_\epsilon(M)$, there there exists
$(\cF,q,J,\rho)$ together with two lagrangian subsheaves $\cL_0$, $\cL_1$
such that 
$$[\omega]=[\tilde \ch(\nabla^{\oplus_0},\nabla^{F_\C})-\tilde
\ch(\nabla^{\oplus_1},\nabla^{F_\C})]=[\tilde \ch(\nabla^{\oplus_0},\nabla^{\oplus_1})]$$ in $H^{odd}(M)$,
where $\nabla^{\oplus_i}$, $i=0,1$ are defined using the decompositions
$F=L_i\oplus J(L_i)$. 
 The right-hand side belongs to $\ch(K^{-1}(X))$ and is therefore rational. On the other hand it depends continuously on $J$. Therefore it is independent of $J$. First reducing to the case that $L_0\cap L_1=\{0\}$ and then choosing
$J$ such that $JL_i=L_{1-i}$, $i=0,1$, we see that $[\omega]=0$.

We now construct a natural morphism of $L(M)$-modules $\bar \gamma :\bar
L(M)\rightarrow K^{-1}_{\R/Z}(M)$. Here we use the definition of
$K^{-1}_{\R/Z}(M)$ in terms of generators and relations given in
\cite{lott9?}, Def. 5 + 6. A generator of $K^{-1}_{\R/Z}(M)$ is a tuple
$(E,h^E,\nabla^E,\rho)$, where $E$ is a $\Z_2$-graded complex vector bundle of
vertual dimension zero, $h^E$ is a hermitean metric and $\nabla^E$ is a metric
connection, both being compatible with the grading, and $\rho\in
\Omega^{odd}(M)/\im(d)$ satisfies $d\rho=\ch(\nabla^E)$. The relations of 
$K^{-1}_{\R/Z}(M)$ are generated by \begin{enumerate}
\item {\bf isomorphism}
$(E,h^E,\nabla^E,\rho)\sim
(E^\prime,h^{E^\prime},\nabla^{E^\prime},\rho)$
if there exists an isomorphism from $E$ to $E^\prime$ which is compatible
with metrics and connections.
\item {\bf direct sum}
$(E,h^E,\nabla^E,\rho)+
(E^\prime,h^{E^\prime},\nabla^{E^\prime},\rho^\prime)=(E\oplus
E^\prime,h^{E\oplus E^\prime},\nabla^{E\oplus E^\prime},\rho+\rho^\prime)$
\item {\bf change of connections}
$(E,h^E,\nabla,\rho)\sim (E,h^E,\nabla^\prime,\rho^\prime)$ if
$\rho^\prime=\rho+\tilde \ch(\nabla^\prime,\nabla)$.
\item {\bf trivial elements}
If $(E,h^E,\nabla^E)$ is a $Z_2$-graded hermitean vector bundle with
connection, then
$(E\oplus E^{op}, h^{E\oplus E},\nabla^{E\oplus E},0)\sim 0$,
where $E^{op}$ denotes $E$ with the opposite grading.
\end{enumerate}
Let $[E,h^E,\nabla^E,\rho]$ denote the class of $(E,h^E,\nabla^E,\rho)$ in
$K^{-1}_{\R/Z}(M)$.

We define $\bar L(M) \rightarrow K^{-1}_{\R/Z}(M)$ by
$\bar\gamma[\cF,q,J,\rho]=[F_\C,h^{F_\C},\nabla^{F_\C,J},\rho]$, where
$h^{F_\C}$ is the hermitean extension of the metric $J^*\circ Q$ on $F$.

The following diagram commutes
$$\begin{array}{ccc}
\bar L(M)&\rightarrow &L(M)\\
\downarrow&&\downarrow\\
K^{-1}_{\R/Z}(M)&\rightarrow & K^0(M)
\end{array}\ .$$
We therefore obtain a natural map
$\bar L(M)\rightarrow \bar L^{\R/\Z}(M)$ which is
in fact surjective. 

We now extend the functor $\bar L$ to $\cTop_{smooth}$  by setting
$$\bar L(X):=\lim_{\stackrel{\longleftarrow}{f:M\rightarrow X}}\bar L(M)\ ,$$
where the limit is taken over the category of manifolds over $X$.
This extension has all functorial properties discussed above.

\section{Eta invariants and index maps}
 
\subsection{The $\eta$-invariant}

Using the $\eta$-invariant of the twisted signature operator for a closed
odd-dimensional oriented manifold $M$ in \cite{lott99} Lott constructed a group
homomorphism $\eta^{Lott}:\bar L^{Lott}(M)\rightarrow \R$.
Its reduction modulo $\Z$ factors over the homomorphism
$\eta^{\R/\Z}:K^{-1}_{\R/\Z}(M)\rightarrow \R/\Z$ constructed in
\cite{lott9?}, which is given by the pairing with the $K$-homology class
induced by the odd signature operator.

Unfortunately, the homomorphism $\eta^{Lott}$ does not factor over the quotient
$\bar L^{Lott}(M)\rightarrow \bar L(M)$. In order to fix this for a given
closed odd-dimensional oriented manifold $M$ we define an extension
$$0\rightarrow \Z\rightarrow \bar L_\epsilon^{ex}(M)\rightarrow
L_\epsilon(M)\rightarrow 0\ ,$$ such that $$\bar\eta:\bar
L_\epsilon^{ex}(M)\rightarrow \bar L_\epsilon(M)\rightarrow
K^{-1}_{\R/\Z}(M)\stackrel{\eta^{\R/\Z}}{\rightarrow} \R/\Z$$ lifts to
$\eta:\bar L_\epsilon^{ex}(M)\rightarrow \R$. Note that $\bar L^{ex}(M)$ is
not a functor on $M$.
Here $n:=\dim(M)$,  $\epsilon_n:=[\frac{n(n-1)}{2}]\in\Z_2$,
and we assume that $\epsilon=1-\epsilon_n$.

Let $(\cF,q)$ be a generator of $L_\epsilon(M)$ and $\cL$ be a lagrangian
locally constant subsheaf of $\cF$. We define an integer $\tau(\cF,q,\cL)$
by the following construction. 
We consider the complex of sheaves
$$\cK^\bullet:\cL\stackrel{i}{\rightarrow} \cF\stackrel{i^*\circ
q}{\rightarrow} \cL^*\ .$$ Let $(E^{\bullet\bullet}_r,d_r)$, $r\ge 1$  be the
associated hyper cohomology spectral sequence. The duality $q$ and the pairing
between $\cL$ and $\cL^*$ induce a duality $q_{\cK}:
\cK\stackrel{\sim}{\rightarrow} \cK^*[-2]$
(the argument "$[-2]$" indicates that $q_{\cK}$ is a map of degree $-2$). The orientation of $M$ and $q$
induce a duality $q_{E_r}:E_r\stackrel{\sim}{\rightarrow}
E_r^*[-2,-\dim(M)]$ (i.e. of bidegree $(-2,-\dim(M))$). Let $N_r$ denote the
total grading on $E_r$. Then we define the symmetric form $Q_r$ on $E_r$ by
$Q_r(v,w)=q_{E_r}((-1)^{\frac{N_r(N_r-1)}{2}}v)(d_r w)$. 
The integer $\tau(\cF,q,\cL)$ is now given by
$$\tau(\cF,q,\cL):=2(\sign(Q_1) + \sign(Q_2)) \ .$$

We first define an extension $L_{\epsilon}^{ex}(M)$ of $L_\epsilon(M)$.
A generator of the semigroup group $\hat L_\epsilon^{ex}(M)$ is a triple
$(\cF,q,z)$ consisting of a locally constant sheaf of finite-dimensional
$\R$-modules, an $\epsilon$-symmetric duality structure $q$, and an integer
$z$. The semigroup operation is given by 
$$(\cF,q,z)+(\cF^\prime,q^\prime,z^\prime):=(\cF\oplus \cF^\prime,q\oplus
q^\prime,z+z^\prime)\ .$$ Then we form
$L^{ex}_\epsilon(M)$ by introducing the relation 
{\bf lagrangian reduction}
 $$(\cF,q,z)\sim (0,0,z+\tau(\cF,q,\cL))$$
if $\cF$ admits a locally constant lagrangian subsheaf $\cL$.
 
There is an exact sequence
$$0\rightarrow \Z \rightarrow L_\epsilon^{ex}(M)\rightarrow L_\epsilon(M)\rightarrow 0\ ,$$
where the maps are the obvious inclusion and projection.
We now define $\bar L_\epsilon^{ex}(M)$ by the following pull-back diagram:
$$\begin{array}{ccc}
\bar L_\epsilon^{ex}(M) &\rightarrow& L_\epsilon^{ex}(M)\\
\downarrow&&\downarrow\\
\bar L_\epsilon(M) &\rightarrow & L_\epsilon(M)
\end{array}\ .
$$
 Then we have exact sequences
$$0\rightarrow \Z \rightarrow \bar L_\epsilon^{ex}(M)\rightarrow \bar
L_\epsilon(M)\rightarrow 0$$ and
$$H^{4*-(-1)^\epsilon}(M,\R)\rightarrow \bar L_\epsilon^{ex}(M) \rightarrow
L_\epsilon^{ex}(M)\rightarrow H^{4*+1-(-1)^\epsilon}(M,\R)$$
in a natural way.

An element of $\bar L^{ex}_\epsilon(M)$ can be represented by a tuple
$(\cF,q,J^F,\rho,z)$. We obtain a class $[\cF,q,J^F,\rho]\in\bar
L_\epsilon^{Lott}(M)$.  If we set
$$\eta(\cF,q,J^F,\rho,z):=\eta^{Lott}([\cF,q,J^F,\rho])-z\ ,$$ then we have
the following result. \begin{theorem}
The map $(\cF,q,J^F,\rho,z)\mapsto \eta(\cF,q,J^F,\rho,z)$ induces a
well-defined homomorphism $\eta:\bar L_\epsilon^{ex}(M)\rightarrow \R$
which  lifts 
$\bar\eta$.
\end{theorem}
We must show that $\eta$ is well-defined
with respect to lagrangian reduction.
The idea is to consider the $\eta$-invariant of the signature operator
on $M$ twisted with the complex $K=\bundle(\cK)$ and to rescale its
differential.   The corresponding adiabatic limits can be understood by the
methods developed in \cite{bismutcheeger89} and \cite{dai91} without
investing essentially new ideas.

\subsection{The secondary index map}

Let $M\rightarrow B$ be a smooth locally trivial fibre bundle with
even-dimensional closed fibres over a compact base $B$ such that the vertical
bundle $TM/B$ is oriented. We set $n:=\dim(TM/B)$.  

In \cite{lott99} Lott constructed a secondary index map
$\pi_*^{\bar L, Lott}:\bar L_\epsilon^{Lott}(M)\rightarrow
L^{Lott}_{\epsilon+\epsilon_n}(B)$ which fits into the commutative
diagram
$$\begin{array}{ccccc}H^{4*-(-1)^\epsilon}(M,\R)&\rightarrow &\bar
L_\epsilon^{Lott}(M) &\rightarrow& \bar L_\epsilon^{\R/\Z}(M)\\
\pi_*^{sign}\downarrow &&\pi_*^{\bar L,Lott} \downarrow &&\pi^{\bar 
L,\R/\Z}_*\downarrow\\
H^{4*-(-1)^{\epsilon+\epsilon_n}}(M,\R)&\rightarrow
&\bar L^{Lott}_{\epsilon+\epsilon_n}(B)&\rightarrow &
\bar L^{\R/\Z}_{\epsilon+\epsilon_n}(B) \end{array}\ .$$
\begin{theorem} By passing to quotients the map $\pi_*^{\bar L, Lott}$
induces a well-defined secondary index map $$\pi^{\bar L}_*:\bar
L_\epsilon(M)\rightarrow \bar L_{\epsilon+\epsilon_n}(B)\ .$$
\end{theorem}
The idea of the proof consists in investigating the adiabatic limits of the
$\eta$-form of the fibrewise signature operator twisted with the complex
$K$ under rescaling the differential. The arguments are similar to the
case of analytic torsion forms. It suffices to adapt the methods developed in
\cite{bismutlebeau91}, \cite{berthomieubismut94}, \cite{bismut97},
\cite{ma99}, \cite{ma00}, \cite{ma001}, \cite{ma002}.

The following diagram commutes 
$$\begin{array}{ccccc}H^{4*-(-1)^\epsilon}(M,\R)&\rightarrow &\bar
L_\epsilon(M) &\rightarrow& L_\epsilon(M)\\
\pi_*^{sign}\downarrow &&\pi_*^{\bar L} \downarrow &&\pi^{L}_*\downarrow\\
H^{4*-(-1)^{\epsilon+\epsilon_n}}(M,\R)&\rightarrow
&\bar L_{\epsilon+\epsilon_n}(B)&\rightarrow &
L_{\epsilon+\epsilon_n}(B) \end{array}\ .$$
The secondary index map is natural with respect to pull-back of fibre bundles.
Furthermore, if $\pi_1:X\rightarrow X_1$ and $\pi_2:X_2\rightarrow B$ is an
iterated bundle with oriented vertical bundles, and $\pi=\pi_2\circ \pi_1$,
then 
\begin{theorem} $\pi^{\bar L}_*=(\pi^{\bar L}_2)_*\circ (\pi^{\bar L}_1)_*$
\end{theorem}
This follows from the functoriality of the  $\eta$ -form studied in
\cite{ma00} in a similar way as the functoriality of the secondary index
in \cite{bunke00} was deduced from the functoriality of the higher analytic
torsion form \cite{ma99}.

\subsection{The index map for $L^{ex}$ and $\bar L^{ex}$}

Let $\pi:M\rightarrow B$ as above and assume that $B$ is oriented, closed, of
dimension $m$.  Furthermore, assume that $\epsilon=1-\epsilon_n-\epsilon_m$.
Let $(\cF,q)$ be a generator of $L_\epsilon(M)$. We define an integer
$\tau(\cF,q,M\stackrel{\pi}{\rightarrow} B)$ by the following construction.
Let $(E^{\bullet\bullet},d_r)$ be the Leray-Serre spectral sequence converging
to sheaf cohomology $H^*(M,\cF)$ with $E_2^{p,q}=H^p(B,R^q\pi_*\cF)$. The
orientations of $TM/B$ and  $B$ and the duality $q$ induce  dualities
$q_{E_r}:E_r\rightarrow E_r^*[-m,-n]$ which give the Poincar\' e duality on
the limit. We define the symmetric form $Q_r$ on $E_r$ by
$Q_r(v,w)=q_{E_r}((-1)^{\frac{N_r(N_r-1)}{2}}v)(d_rw)$,
where $N_r$ denotes the total degree.
Then we set  
$$\tau(\cF,q,M\stackrel{\pi}{\rightarrow} B):=2\sum_{r\ge 2} \sign(Q_r)\ .$$
\begin{theorem}\label{rr}
\begin{enumerate}
\item
The prescription
$$\pi^{L^{ex}}_*[\cF,q,z]:=[H^\bullet
R\pi_*\cF,\pi_*(q),z-\tau(\cF,q,M\stackrel{\pi}{\rightarrow} B)]\ .$$ defines
an extended index map $\pi^{L^{ex}}_*:L^{ex}_\epsilon(M)\rightarrow
L^{ex}_{\epsilon+\epsilon_n}(B)$. \item The extented index map is functorial
with respect to iterated fibre bundles.
\end{enumerate}
\end{theorem}
The following diagram  
$$\begin{array}{ccccc} \bar
L_\epsilon(M)&\rightarrow&L_\epsilon(M)&\leftarrow&L^{ex}_\epsilon(M)\\
\pi^{\bar L}_*\downarrow &&\pi^L_*\downarrow&&\pi^{L^{ex}}_*\downarrow\\
\bar L_{\epsilon+\epsilon_n}(B)&\rightarrow&
L_{\epsilon+\epsilon_n}(B)&\leftarrow&L^{ex}_{\epsilon+\epsilon_n}(B)
\end{array}$$
commutes and induces an extended secondary index map
$\pi^{\bar L^{ex}}_*:\bar L^{ex}(M)_\epsilon \rightarrow \bar
L^{ex}_{\epsilon+\epsilon_n}(B)$ which is functorial with respect to iterated
fibre bundles. 
We have the following compatibility of this extended secondary index map with
the $\eta$-homomorphism.
\begin{theorem}\label{rr1}
$$\begin{array}{ccc}
\bar L^{ex}_\epsilon(M)&\stackrel{\eta}{\rightarrow}&\R\\
\pi^{\bar L^{ex}}_*\downarrow&&\|\\
\bar L^{ex}_{\epsilon+\epsilon_n}(B)&\stackrel{\eta}{\rightarrow}&\R
\end{array}$$
\end{theorem}
This theorem essentially follows from Dai's adabatic limit formula for
the $\eta$-invariant \cite{dai91}.

\bibliographystyle{plain}

\end{document}